\documentclass[12pt]{article}%
\usepackage{amsmath}
\usepackage{amsfonts}
\usepackage{amssymb}
\usepackage{graphicx}
\usepackage{graphics}%

\def\alp{\alpha}
\def\bta{\beta}
\def\dta{\delta}
\def\upa{\uparrow}

\def\tht{\theta}
\def\ol{\overline}

\def\gam{\gamma}
\def\Gam{\Gamma}

\def\ome{\omega}
\def\tms{\times}
\def\lds{\ldots}

\title{On a Small Elliptic Perturbation of a Backward-Forward Parabolic Problem, with Applications to Stochastic Models}

\author{{\normalsize by}\\ \\
\normalsize Diego Dominici and Charles Knessl\\ \\ \normalsize Dept. of Mathematics, Statistics\\ \normalsize and Computer Science (M/C 249)\\ \normalsize University of Illinois at Chicago\\ \normalsize 851 South Morgan Street\\
\normalsize Chicago, IL 60607-7045}

\date{}

\begin{document}
\maketitle

\begin{abstract}
We consider an elliptic PDE in two variables.  As one parameter approaches zero, this PDE collapses to a parabolic one, that is forward parabolic in a part of the domain and backward parabolic in the remainder.  Such problems arise naturally in various stochastic models, such as fluid models for data-handling systems and Markov-modulated queues.  We employ singular perturbation methods to study the problem for small values of the parameter.
\end{abstract}

\section{Introduction}
\label{sec1}

We consider the following boundary value problem for $F=F(x,\xi)$
\renewcommand{\theequation}{\arabic{section}.\arabic{equation}}
\setcounter{section}{1}
\setcounter{equation}{0}
\begin{eqnarray}
& & DF_{xx}+(a-\xi)F_x+F_{\xi\xi}+(\xi F)_\xi=0;\qquad x>0, \quad -\infty<\xi<\infty\label{eq1.1}\\
& & \quad DF_x(0,\xi)+(a-\xi)F(0,\xi)=0,\quad -\infty<\xi<\infty\label{eq1.2}\\
& & \qquad \int\limits^\infty_{-\infty} \int\limits^\infty_0 F(x,\xi)dx\,d\xi=1.\label{eq1.3}
\end{eqnarray}
Here $a$ and $D$ are positive constants, and $F$ is a bivariate probability density.  This problem arose as an asymptotic limit of a Markov-modulated queueing model, where the input process is generated by $N$ sources that turn ``on'' and ``off'' at exponential waiting times.  When on, a source generates a Poisson arrival stream to the queue.  The joint distribution of the number of on sources and the queue length satisfies a complicated system of difference equations, that as $N\to\infty$ may be approximated by the problem (\ref{eq1.1})-(\ref{eq1.3}).  The variable $x$ is related to the queue length and $\xi$ corresponds to a scaled, centered measure of the number of on sources.  An empty queue has $x=0$, and $\xi=0$ means that the number of active sources equals its mean value.  

The diffusion coefficient $D$ in (\ref{eq1.1}) measures variability effects in the service time distribution, while $a$ measures the difference between average output and input rates to the queue.  The condition $a>0$ guarantees stability and the existence of a steady state distribution.  
Note that $D\to 0$ means that the queue becomes a deterministic, or ``fluid'', process.  
A detailed derivation of (\ref{eq1.1})-(\ref{eq1.3}) can be found in \cite{1:Kn-Ti}.  There we also showed that the solution can be reduced to either a Fredholm integral equation of the second kind (for $F(0,\xi)$), or to finding a single eigenvector of an infinite matrix, whose elements are expressed in terms of Laguerre functions.

The complexity of the exact solution to (\ref{eq1.1})-(\ref{eq1.3}) suggests that an asymptotic analysis may be fruitful.  In \cite{1:Kn-Ti} we considered the limit $D\to\infty$ and discussed numerical methods that work well for $D$ large and $D=O(1)$.  The purpose of this note is to analyze the opposite limit, namely $D\to 0^+$.  This is a very singular limit, as we show in what follows.  Also, the numerical methods involve truncating an infinite matrix, say to an $M\tms M$ matrix. For $D\to 0^+$ the eigenvector we seek decays on the scale (see [1, section 4.2]) $O(D^{-2})$.  Thus for say $D=.1$, $M$ must be of the order of about 500 before an accurate result is obtained.  The basic matrix is not sparse; it does become diagonally dominant for $D\to\infty$, but not for $D\to 0^+$.

Let us briefly discuss the problem with $D=0$.  Denoting by ${\cal F}(x,\xi)$ the solution to (\ref{eq1.1}) with $D=0$, we see that the elliptic PDE degenerates into a parabolic one, that is forward parabolic in the range $\xi>a$ and backward parabolic for $\xi<a$.  The boundary condition (\ref{eq1.2}) becomes ${\cal F}(0,\xi)=0$ and can only be applied in the range $\xi>a$, where the PDE is forward parabolic.  The study of backward/forward parabolic problems dates as far back as 1914 (see Gerrey \cite{2:Ge}) and particular problems of the type here are analyzed in [3,4,5].  In \cite{4:Kn-Mo} we give an explicit expression for ${\cal F}(x,\xi)$, subject to the condition ${\cal F}(\infty,\xi)=(2\pi)^{-1/2}e^{-\xi^2/2}$ (replacing (\ref{eq1.3})).  The function ${\cal F}$ can also be interpreted as an asymptotic limit of a probability distribution for a discrete stochastic model, that was formulated and analyzed in \cite{6:An-Mi-So}.  Here we are interested in how the solution of the elliptic PDE approaches that of the parabolic one.  For the latter the values of ${\cal F}(0,\xi)$ for $\xi<a$ are unknown and must be computed as a part of the solution.  Probabilistically, this corresponds to boundary mass along this part of the boundary.  For any $D>0$, (\ref{eq1.3}) shows that $F$ is a proper density function.  Now the ``no-flux'' boundary condition (\ref{eq1.2}) applies for all $\xi$, and no boundary mass develops.  We also note that by integrating (\ref{eq1.1}) over $x\in(0,\infty)$ and using (\ref{eq1.2}), the marginal distribution $\int^\infty_0F(x,\xi)dx$ satisfies an elementary ODE that is readily solved to give
\begin{equation}
\int\limits^\infty_0 F(x,\xi)dx=\frac{1}{\sqrt{2\pi}}e^{-\xi^2/2}.
\label{eq1.4}
\end{equation}
Here we also used the normalization (\ref{eq1.3}).  Since $F$ is a density in $x$ while ${\cal F}$ is a distribution in $x$, we expect to compare $F$ to ${\cal F}_x$ for $x>0$.

To get a qualitative picture of the solution, we plot in Figure 1 the deterministic approximation to (\ref{eq1.1}).  

\begin{figure}[t]
\begin{center}
\rotatebox{0} {\resizebox{15cm}{!}{\includegraphics{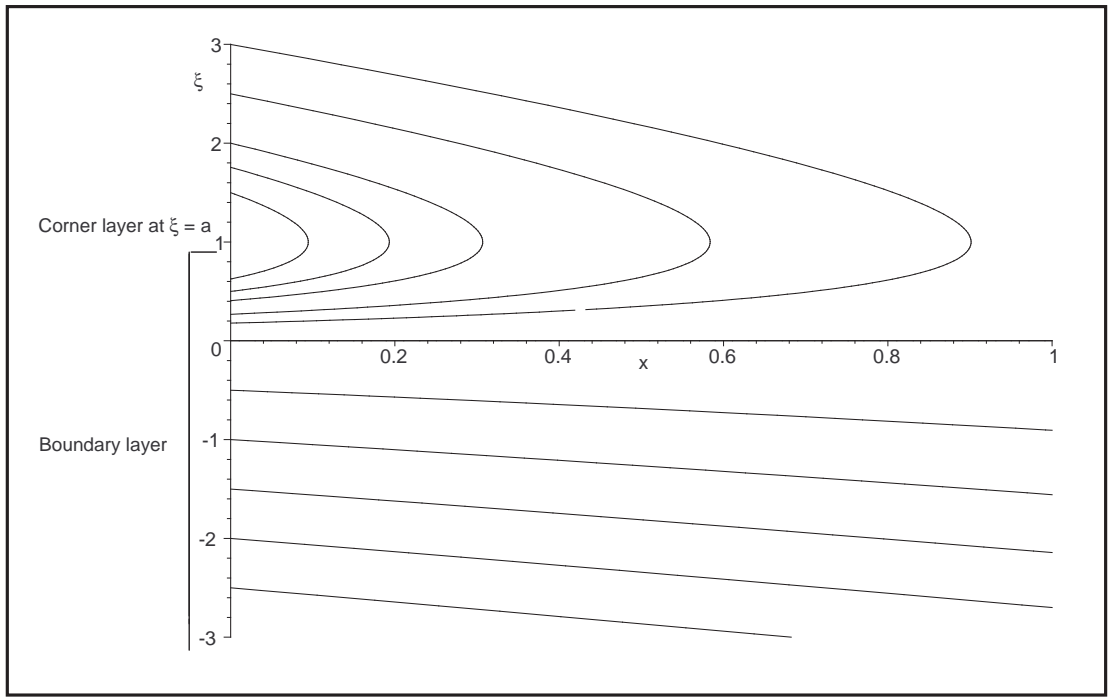}}}
\end{center}
\caption{A sketch of the deterministic trajectories for $a = 1$.}%
\end{figure}

This corresponds to neglecting diffusion in both the $x$- and $\xi$-variables, and is the phase-plane flow $\dot x=\xi-a$, $\dot\xi=-\xi$.  We plot a few trajectories for $x>0$, and note that once they hit $x=0$ (necessarily with $\xi<a$), they stay at $x=0$ and flow to the equilibrium at $\xi=0$.  The figure also indicates the boundary and corner layers that arise in the analysis for $D\to 0$.  

\section{Singular Perturbation Analysis}
\label{sec2}

We shall analyze (\ref{eq1.1})-(\ref{eq1.3}) for $D\to0^+$ in the three ranges (i) $x>0$, $-\infty<\xi<\infty$ (or $x\approx 0$, $\xi>a$); (ii) $x\approx 0$, $\xi<a$; and (iii) $x\approx 0$, $\xi\approx a$.

\subsection{Outer Solution}
\label{subsec2.1}

We assume $x>0$ and expand $F(x,\xi)$ as $F(x,\xi)=F_0(x,\xi)+DF_1(x,\xi)+O(D^2)$ to find that
\setcounter{section}{2}
\setcounter{equation}{0}
\begin{equation}
(\xi-a)F_{0,x}=F_{0,\xi\xi}+(\xi F_0)_\xi.
\label{eq2.1}
\end{equation}
We note that when $\xi>a$, (\ref{eq2.1}) is a diffusion equation with $x$ taking the place of time, and when $\xi<a$ it is a diffusion equation run in reverse time.  As is well-known [2,3,5], we need to impose boundary conditions for (\ref{eq2.1}) on the ``incoming'' half of the boundary, where $x=0$ and $\xi>a$.  We can derive the appropriate ``half boundary condition'' for (\ref{eq2.1}) by considering a boundary layer about $x=0$.

We thus consider the boundary layer where $x=O(D)$.  Setting $x=Du$ and $F(x,\xi)=G(u,\xi)\sim D^{\nu_0}G_0(u,\xi)$, we find that
\begin{eqnarray}
& & G_{0,uu}+(a-\xi)G_{0,u}=0; \ \ u>0, \qquad -\infty<\xi<\infty\label{eq2.2}\\
& & G_{0,u}(0,\xi)+(a-\xi)G_0(0,\xi)=0, \,\qquad -\infty<\xi<\infty\label{eq2.3}
\end{eqnarray}
and hence $G_0$ is proportional to $\exp\bigl[(\xi-a)u\bigr]$.  For $\xi>a$ this grows exponentially as $u\to\infty$ and cannot match to the outer solution $F_0(x,\xi)$ as $x\to 0^+$, for any choice of $\nu_0$.  We thus conclude that no boundary layer develops along $\xi>a$, where (\ref{eq2.1}) is forward parabolic.  Thus $F_0$ must satisfy $F_0(0,\xi)=0$ for $\xi>a$, which is the limit of (\ref{eq1.2}) as $D\to 0^+$.

The explicit solution to (\ref{eq2.1}), subject to the above ``half boundary condition,'' is given in \cite{4:Kn-Mo}.  
There we solved (\ref{eq2.1}) using a Laplace transform in the $x$-variable.  Due to the forward/backward nature of the PDE, the spectrum has both positive and negative eigenvalues.  Insisting that the solution decay as $x\to\infty$, the positive eigenvalues must be absent.  This, together with the half boundary condition for $\xi>a$ is sufficient to determine the solution. 

Noting that $F_0$ should correspond to ${\cal F}_x$, we give two different representations of this solution.  The first is the contour integral
\begin{eqnarray}
\qquad 
F_0(x,\xi) & = & K\frac{1}{2\pi i}\int_{\rm Br}\frac{e^{-\xi^2/4}}{\sqrt{2\pi}} D_{a\tht+\tht^2}(\xi+2\tht)A(\tht)e^{\tht x}d\tht,\label{eq2.4}\\
A(\tht) & = & \frac{a}{\tht+a}\exp\left[-\tht\zeta\bigl(\frac12\bigr)+\frac12\gam\tht(\tht+a)\right]\prod\limits^\infty_{m=1}\frac{\dta_m}{\tht+\dta_m}\exp\left[\frac{\tht}{\sqrt m}-\frac{\tht(\tht+a)}{2m}\right],\nonumber\\
\dta_n & = & \frac12\bigl[a+\sqrt{a^2+4n}\,\bigr], \quad n=0,1,2,\lds \ .\nonumber
\end{eqnarray}
Here $D_p(\cdot)$ is the parabolic cylinder function of order $p$, $\gam$ is Euler's constant and $\zeta (\cdot)$ is the Riemann zeta function.  The Bromwich contour Br in (\ref{eq2.4}) is any vertical contour in the half-plane ${\rm Re}(\tht)>-a$, where the integrand is analytic. 
Note that the only singularities of the integrand are simple poles at $\tht=-\dta_n$, $n\ge 0$.  If we close the contour to the left, then we can convert the integral to the sum of the residues at the poles.  This leads to a second representation of the solution in terms of the problem's eigenfunctions:

\begin{eqnarray}
\qquad
F_0(x,\xi) & = & K\sum\limits^\infty_{n=0}c_ne^{-\dta_n x}\frac{e^{-\xi^2/2}}{\sqrt{2\pi}}e^{\xi\dta_n}H_n\left(\frac{\xi-2\dta_n}{\sqrt 2}\right),\label{eq2.5}\\
c_n & = & \frac{a2^{-n/2}e^{-\dta^2_n}}{a-\dta_n}\exp\left[\dta_n\zeta \bigl(\frac12\bigr)+\frac12 n\gam\right]{\prod^{\infty}_{m=1}}' \frac{\dta_m}{\dta_m-\dta_n}\exp\left[-\frac{\dta_n}{\sqrt m}-\frac{n}{2m}\right],\nonumber\\
& &  \qquad n\ge 1,\nonumber\\
c_0 & = & ae^{-a^2}\exp\left[a\zeta\bigl(\frac12\bigr)\right]\prod\limits^\infty_{m=1} \frac{\dta_m}{\dta_m-a}\exp\left(-\frac{a}{\sqrt m}\,\right).\nonumber
\end{eqnarray}
The product in the expression for $c_n$ omits the term $m=n$ and $H_n(\cdot)$ is the $n^{\rm th}$ Hermite polynomial.  

The constant $K$ is a normalization constant that will ultimately be determined from (\ref{eq1.3}).  We cannot determine it at this stage since we shall show that there is also an $O(1)$ probability mass in a boundary layer near $x=0$ (with $\xi<a$), that contributes to (\ref{eq1.3}).  

\subsection{Boundary Layer}
\label{subsec2.2}

While no boundary layer develops for $\xi>a$, there is one for $\xi-a<0$ and $x=O(D)$. We again scale $x=Du$, and set
\begin{equation}
F(x,\xi)=G(u,\xi)=\frac{1}{D} G_0(u,\xi)+G_1(u,\xi)+O(D).
\label{eq2.6}
\end{equation}
Then $G_0$ will satisfy the problem (\ref{eq2.2}) and (\ref{eq2.3}), whose general solution is 
\begin{equation}
G_0(u,\xi)=g_0(\xi)e^{(\xi-a)u},
\label{eq2.7}
\end{equation}
and this decays exponentially as $u\to\infty$ for $\xi<a$.  We will determine $g_0(\xi)$ shortly, and also show that the expansion (\ref{eq2.6}) must include the term of order $O(D^{-1})$ in order to match to the outer expansion.

Using (\ref{eq2.6}) in (\ref{eq1.1}) we find that the correction term $G_1$ satisfies
\begin{eqnarray}
& & G_{1,uu}+(a-\xi)G_{1,u}=-G_{0,\xi\xi}-(\xi G_0)_\xi\label{eq2.8}\\
& & \qquad =-e^{(\xi-a)u}\bigl[g''_0+(\xi+2u)g'_0+(u^2+u\xi+1)g_0\bigr]\nonumber
\end{eqnarray}
and (\ref{eq1.2}) yields
\begin{equation}
G_{1,u}(0,\xi)+(a-\xi)G_1(0,\xi)=0.
\label{eq2.9}
\end{equation}
We write the solution to (\ref{eq2.8}) as
\begin{equation}
G_1(u,\xi)=e^{(\xi-a)u}\ol G(u,\xi)+h(\xi).
\label{eq2.10}
\end{equation}
Then the BC (\ref{eq2.9}) yields
\begin{equation}
\ol G_u(0,\xi)+(a-\xi)h(\xi)=0.
\label{eq2.11}
\end{equation}
From (\ref{eq2.8}) and (\ref{eq2.10}) we obtain
$$\ol G_{uu}+(\xi-a)\ol G_u=-\bigl[g''_0+(\xi g_0)'+u(2g'_0+\xi g_0)+u^2g_0\bigr]$$
so that
\begin{equation}
\ol G(u,\xi)=\alp(\xi)\frac{u^3}{6}+\bta(\xi)\frac{u^2}{2}+\gam(\xi)u+\dta(\xi)
\label{eq2.12}
\end{equation}
where
\begin{eqnarray}
\alp(\xi) & = & -\frac{2}{\xi-a}g_0(\xi),\quad \bta(\xi)=-\frac{2}{\xi-a}g'_0(\xi)+\left[\frac{2}{(\xi-a)^2}-\frac{\xi}{\xi-a}\right]g_0(\xi),\label{eq2.13}\\
\gam(\xi) & = & -\frac{1}{\xi-a}\bigl[g''_0+(\xi g_0)'\bigr]+\frac{1}{(\xi-a)^2}[2g'_0+\xi g_0]-\frac{2}{(\xi-a)^3}g_0\nonumber
\end{eqnarray}
and $\dta(\xi)$ is an arbitrary function.

The asymptotic matching of the inner and outer solutions requires that as $u\to\infty$ (\ref{eq2.6}) agrees with $F_0$ as $x\to 0^+$, and since $\xi<a$ this implies that
\begin{equation}
h(\xi) =F_0(0,\xi)=K\frac{1}{2\pi i}\int\limits_{\rm Br}\frac{e^{-\xi^2/4}}{\sqrt{2\pi}} D_{a\tht+\tht^2}(\xi+2\tht)A(\tht)d\tht.
\label{eq2.14}
\end{equation}
Since $\ol G_u(0,\xi)=\gam(\xi)$, (\ref{eq2.11}) yields
\begin{equation}
\gam(\xi)=(\xi-a)h(\xi)=\frac{d^2}{d\xi^2}\left(\frac{g_0}{a-\xi}\right)+\frac{d}{d\xi}\left(\frac{\xi g_0}{a-\xi}\right)
\label{eq2.15}
\end{equation}
where the last equality follows from (\ref{eq2.13}).  We thus set $g_0(\xi)=(a-\xi) \ol g(\xi)$ to find that
\begin{equation}
\ol g''(\xi)+(\xi \ol g)'(\xi)+(a-\xi)F_0(0,\xi)=0.
\label{eq2.16}
\end{equation}
Also, the leading term in the boundary layer becomes
\begin{equation}
G(u,\xi)\sim \frac{1}{D} (a-\xi)e^{(\xi-a)u}\ol g(\xi), \quad \xi<a.
\label{eq2.17}
\end{equation}
By integrating (\ref{eq2.1}) from $x=0$ to $x=\infty$ and comparing the result to (\ref{eq2.15}), we conclude that ${\cal D}(\xi)\equiv \int^\infty_0 F_0(x,\xi)dx+\ol g(\xi)$ satisfies ${\cal D}''(\xi)+(\xi{\cal D})'(\xi)=0$ for $\xi <a$, so we have
\begin{equation}
\ol g(\xi)+\int^\infty_0F_0(x,\xi)dx=\frac{L}{\sqrt{2\pi}}e^{-\xi^2/2}, \quad \xi<a
\label{eq2.18}
\end{equation}
for some constant $L$.  But the left side of (\ref{eq2.18}) is the leading term in the expansion of the left side of (\ref{eq1.4}) for $D\to 0^+$, hence $L=1$.  Note also that the normalization condition in (\ref{eq1.3}) implies to leading order that
\begin{equation}
\int\limits^a_{-\infty} \ol g(\xi)d\xi+\int\limits^\infty_{-\infty}\int\limits^\infty_0 F_0(x,\xi)dx\,d\xi=1.
\label{eq2.19}
\end{equation}
It remains only to determine $K$.  For $\xi>a$ there is no boundary layer contribution for $x=O(D)$ to (\ref{eq1.4}) and thus 
\begin{equation}
\int\limits^\infty_0 F_0(x,\xi)dx=K\frac{1}{2\pi i} \int_{{\rm Br}'} \frac{e^{-\xi^2/4}}{\sqrt{2\pi}} D_{a\tht+\tht^2}(\xi+2\tht)\frac{A(\tht)}{\tht}d\tht=\frac{1}{\sqrt{2\pi}}e^{-\xi^2/2}
\label{eq2.20}
\end{equation}
where we restrict $-a< {\rm Re}(\tht)<0$ on ${\rm Br}'$ (this allows us to perform the $x$-integration using (\ref{eq2.4}) by simply replacing $e^{\tht x}$ by $1/\tht$).  For $\xi>a$ we can close the integrating contour in the right half-plane, picking up the residue from the pole at $\tht=0$.  We conclude that $K=1$, as $A(0)=1$ and $D_0(\xi)=e^{-\xi^2/4}$.

To summarize, we have obtained the outer solution as (\ref{eq2.4}) with $K=1$ and the boundary layer as (\ref{eq2.17}), with (\ref{eq2.18}) and $L=1$.  The correction term in the boundary layer is given by $G_1$ in (\ref{eq2.6}), with (\ref{eq2.10})-(\ref{eq2.14}).  The function $\dta(\xi)$ can ultimately be found by higher order matchings; this would require we compute the correction term $F_1$ in the outer solution.  The analysis shows that for $D$ small and $x>0$ (or $x\approx 0$ with $\xi>a$) the shape of the density for the elliptic problem is, to leading order, the same as that for the parabolic one.  However for $x$ small and $\xi<a$ a boundary layer develops for the elliptic problem, and this takes the place of the boundary mass present in the parabolic problem.

\subsection{Corner Layer}
\label{subsec2.3}

Form (\ref{eq2.17}) we see that the leading term vanishes as $\xi\upa a$, indicating another non-uniformity in the asymptotics.  We thus study how the boundary layer disappears as $\xi$ increases through the critical value $a$, by introducing the scaling
\begin{equation}
\xi -a=D^{1/4}z,\qquad x=D^{3/4} X
\label{eq2.21}
\end{equation}
with
\begin{equation}
F(x,\xi)=H(X,z) \sim D^\nu H_0(X,z)
\label{eq2.22}
\end{equation}
where $\nu$ is a constant to be determined.  
From (\ref{eq1.1}) and (\ref{eq1.2}) we find that $H_0$ satisfies
\begin{eqnarray}
& & H_{0,XX}-zH_{0,X} +H_{0,zz}=0;\qquad X>0, \quad -\infty<z<\infty\label{eq2.23}\\
& & \qquad H_{0,X} (0,z)-zH_0(0,z)=0,\quad -\infty<z<\infty.\label{eq2.24}
\end{eqnarray}
We note that $u=D^{-1/4}X$ and $(\xi-a)u=zX$.  
Note also that going from the $u$ to $X$ scales represents a ``thickening'' of the boundary layer as $\xi\upa a$.  

We seek solutions of (\ref{eq2.23}) in the separable form $H_0=e^{\tau X}{\cal H}(z;\tau)$, where $\tau\in{\Bbb C}$ is a separation constant. We  find that ${\cal H}$ satisfies the Airy equation
$${\cal H}_{zz}=(z\tau-\tau^2){\cal H},$$
so that ${\cal H}$ will be proportional to $Ai\bigl(\tau^{1/3}(z-\tau)\bigr)$ (we expect $H_0$ to decay as $z\to\infty$, corresponding to the disappearance of the boundary layer).  We thus argue that a general solution to (\ref{eq2.23}) will have the form 
\begin{equation}
H_0(X,z)=\frac{1}{2\pi i}\int\limits_C e^{\tau X}Ai(\tau^{1/3}z-\tau^{4/3})f(\tau)d\tau
\label{eq2.25}
\end{equation}
for some function $f$ and contour $C$ in the $\tau$-plane.  If the contour $C$ can be chosen in the range Re$(\tau)<0$ then we integrate (\ref{eq2.23}) and use (\ref{eq2.24}) and (\ref{eq2.25}) to obtain
\begin{equation}
\frac{d^2}{dz^2}\left\{\frac{1}{2\pi i} \int\limits_C \frac{f(\tau)}{\tau} Ai(\tau^{1/3}z-\tau^{4/3})d\tau\right\}=0,
\label{eq2.26}
\end{equation}
so that the integral must be a linear function of $z$, and $f$ satisfies a Fredholm integral equation of the first kind.

While we have not been able to determine $f(\tau)$ explicitly, we can show that (\ref{eq2.25}) can asymptotically match with the outer solution in (\ref{eq2.4}).  For $X$ and $z$ large and positive we expect to have to scale $\tau$ in (\ref{eq2.25}) to be small.  If $f(\tau)\sim f_1\tau^{\nu_1}$ as $\tau\to 0$ the right side of (\ref{eq2.25}) may be approximated by
\begin{equation}
\frac{1}{2\pi i}\int\limits_C f_1\tau^{\nu_1}e^{\tau X}Ai(\tau^{1/3} z)d\tau.
\label{eq2.27}
\end{equation}

We next examine (\ref{eq2.4}) for $x\to 0$ and $\xi\to a$.  For $x\to 0$ we must consider $\tht\to\infty$.  In \cite{4:Kn-Mo} we evaluated the infinite product in (\ref{eq2.4}) as $\tht\to\infty$, with the result
\begin{equation}
A(\tht)\sim\frac{a}{\tht} e^{C_3}e^{-a^2/4}\sqrt\tht e^{-(\tht^2+a\tht)\log\tht}e^{\tht^2/2}.
\label{eq2.28}
\end{equation}
Here $C_3=C_3(a)$ is a constant.  With the scaling $b\to\infty$, $Y\to\infty$ and  $Y-2\sqrt b=O(b^{-1/6})$, the parabolic cylinder function $D_b(Y)$ may be approximated by an Airy function:
\begin{equation}
D_b(Y)\sim 2^{b/2}\Gam\left(\frac{b+1}{2}\right)b^{1/6}Ai\bigl(b^{1/6}(Y-2\sqrt b\,)\bigr).
\label{eq2.29}
\end{equation}
Applying (\ref{eq2.29}) with $b=a\tht+\tht^2$ and $Y=\xi+2\tht$ we have 
$$(Y-2\sqrt b\,)^{1/6} = [\xi+2\tht-2\sqrt{\tht^2+a\tht}\,][a\tht+\tht^2]^{1/6}\sim (\xi-a)\tht^{1/3}$$
and then
\begin{equation}
D_{a\tht+\tht^2}(\xi+2\tht)\sim\sqrt{2\pi}\tht^{1/3}e^{(1/2)(\tht^2+a\tht)\log(\tht^2+a\tht)}e^{-(1/2)(a\tht+\tht^2)}Ai\bigl((\xi-a)\tht^{1/3}\bigr).
\label{eq2.30}
\end{equation}
This holds in the limit $\tht\to\infty$ with $\xi-a=O(\tht^{-1/3})$.  Combining (\ref{eq2.28}) and (\ref{eq2.30}) in (\ref{eq2.4}) and recalling that $K=1$ leads to
\begin{eqnarray}
F_0(x,\xi) & \sim & \frac{1}{2\pi i}\int\limits_{\rm Br} e^{\tht x} ae^{C_3} e^{-a^2/4}\frac{Ai((\xi-a)\tht^{1/3})}{\tht^{1/6}}d\tht\label{eq2.31}\\
& = & ae^{C_3}e^{-a^2/4}\left[\frac{1}{2\pi i}\int\limits_{\rm Br} e^{\tau X} \frac{Ai(z\tau^{1/3})}{\tau^{1/6}}d\tau\right] D^{-5/8}.\nonumber
\end{eqnarray}
By comparing $D^\nu\tms$ (\ref{eq2.27}) with (\ref{eq2.31}) we see that matching is possible with
\begin{equation}
\nu=-\frac{5}{8}, \ \nu_1 = -\frac16, \ f_1=ae^{C_3} e^{-a^2/4}.
\label{eq2.32}
\end{equation}
Thus the density in this corner region will be $O(D^{-5/8})$.  This means that the total mass in the corner is (noting that $dx\,d\xi=Ddz\,dX$) $O(D^{3/8})$.  

Finally we investigate the matching between the corner and boundary layers.  Noting that $\ol g(\xi)$ represents the mass in the boundary layer for a given $\xi<a$, we expect that $\ol g\to 0$ as $\xi\upa a$.  Let us assume that
$$\ol g(\xi)\sim \ol g_2(a-\xi)^{\nu_2},\qquad \xi\to a.$$
Then since $a-\xi=-D^{1/4}z$ the boundary layer becomes $O\bigl[D^{-1}(a-\xi)^{\nu_2+1}\bigr]=O[D^{(\nu_2-3)/4}]$ on the $z$ scale.  This can match to the corner layer if $(\nu_2-3)\big/4=-5/8$, i.e., $\nu_2=1/2$.  Now suppose we let $z\to-\infty$ and $X\to 0$ in (\ref{eq2.24}), scaling $\tau=z+\ome(-z)^{-1/3}$.  This yields
\begin{eqnarray}
D^{-5/8}H_0(X,z) & = & D^{-5/8} e^{zX}(-z)^{-1/3} \frac{1}{2\pi i}\int\limits_{C'} e^{\ome X/(-z)^{1/3}}f(z+\ome(-z)^{1/3})\label{eq2.33}\\
& & \qquad \tms Ai\bigl\{[-z-\ome(-z)^{-1/3}\bigr]^{1/3}\ome (-z)^{1/3}\bigr\} d\ome.\nonumber
\end{eqnarray}
If $f(z)\sim f_3(-z)^{\nu_3}$ as $z\to-\infty$ the left side of (\ref{eq2.33}) should approach
\begin{equation}
D^{-5/8}e^{zX}f_3(-z)^{\nu_3-1/3}\left[\frac{1}{2\pi i} \int\limits_{C'} Ai(\ome)d\ome\right].
\label{eq2.34}
\end{equation}
This can match to (\ref{eq2.16}) if
$$\nu_3=\frac{11}{6},\qquad \ol g_2=f_3\left[\frac{1}{2\pi i}\int\limits_{C'}Ai(\ome)d\ome\right].$$
This shows that the matchings between the corner layer and other two expansions are possible if $f(\tau)=O(\tau^{-1/6})$ as $\tau\to 0$ and $f(\tau)=O\bigl(|\tau|^{11/6}\bigr)$ as $\tau\to-\infty$.

\section{Conclusion}
\label{sec3}

Our analysis showed that the density $F(x,\xi)$ is $O(1)$ for $x>0$, and is large $(O(D^{-1}))$ in the boundary layer where $x=O(D)$ and $\xi<a$.  It was also shown to be large $(O(D^{-5/8}))$ in the corner layer where $(x,\xi)\approx (0,a)$, with the precise scaling given by (\ref{eq2.21}).  The mass in the boundary layer and outer region are comparable, while that in the corner layer is asymptotically small.  We compare this to the solution of the forward/backward problem in \cite{4:Kn-Mo}, where there is non-zero boundary mass along $x=0$ and $\xi<a$.  

We are presently investigating the problem (\ref{eq1.1})-(\ref{eq1.3}) in the limit $a\to\infty$, with $D$ fixed.  This corresponds, after some rescaling, to having small diffusion in both the $x$ and $\xi$ directions.  Some preliminary results show that a geometrical optics type expansion is fruitful, and that the problem yields interesting asymptotic structures, such as a caustic boundary and interior cusped caustics.

\section*{Acknowledgement}

The work of Charles Knessl was partially supported by NSF grants DMS 99-71656 and DMS 02-02815, and NSA grant MDA 904-03-1-0036.  The work of Diego Dominici was supported by NSF grant DMS 99-73231.  We thank Professor Floyd Hanson for providing the support from the latter grant.

\end{document}